\documentclass[12pt]{amsart}
\usepackage{graphicx,pinlabel}
\newtheorem{theorem}{Theorem}[section]

\newtheorem*{maintheorem1}{Theorem~\ref{thm:main1}}
\newtheorem*{maintheorem2}{Theorem~\ref{thm:main2}}
\newtheorem*{corollary1}{Corollary~\ref{cor:1}}
\newtheorem{corollary}[theorem]{Corollary}
\newtheorem{lemma}[theorem]{Lemma}
\newtheorem{proof of theorem}[theorem]{proof of theorem 1.1}
\newtheorem{definition}[theorem]{Definition}
\newcommand{\Out}{\operatorname{Out}}
\newcommand{\Aut}{\operatorname{Aut}}

\newcommand{\PP}{\operatorname{\mathbb{P}}}
\newcommand{\SP}{\operatorname{\mathbb{S}}}

\begin{document}
\title{Normal Tori in $\sharp_n (S^2\times S^1)$}
\author{Funda G\"{U}LTEPE}
\address{Department of Mathematics\\
University of Oklahoma\\ Norman, OK 73019}
\email{fgultepe@math.ou.edu}
\urladdr{http://www2.math.ou.edu/~fgultepe/}
\begin{abstract}The fundamental group of $M = \sharp_n (S^2\times S^1)$ is  $F_n$, the free group with $n$ generators.
 There is a 1-1 correspondence between the equivalence classes of $\mathbb{Z}$-- splittings of $F_n$ and homotopy classes of embedded essential tori in $M$. We define and prove a local notion of minimal intersection of a torus with respect to a maximal sphere system in $M$, which generalizes Hatcher's work \cite{H1} on 2-spheres in the same manifold.
\end{abstract}
\renewcommand{\t}{\widetilde{t}/\gamma}
\newcommand{\tprime}{\widetilde{t'}/\gamma}

\maketitle
\section{introduction}
The study of the group of outer automorphisms $\Out(F_n)$ of the free group $F_n$ on $n$ letters is closely related to the study of the spaces on which  it acts. One such space is Culler-Vogtmann space, or ``Outer Space''. It was first introduced by Culler and Vogtmann in \cite{Vogt1} and it is based on regarding $F_n$ as the fundamental group of a graph. Another such space is obtained from a 3-manifold, $\sharp_n (S^2\times S^1)$ , the connected sum of $n$ copies of $S^2\times S^1$, which we will denote by $M$ for the rest of the paper. The fundamental group of $M$ is also $F_n$ and $\Out(F_n)$ acts on the sphere complex, which is a simplicial complex whose simplices correspond to systems of 2-spheres in $M$. Hatcher and Vogtmann used the sphere complex to prove a homological stability property of $\Aut(F_n)$ in \cite{HatcCerf}. Also, Gadgil \cite{Gad1} gave an algorithm to determine the embedded spheres in $M$ which correspond to the splittings of the free group $F_n$. In \cite{Gad2} Gadgil and Pandit stated  and proved a geometric and algebraic intersection number argument for the sphere complex of $M$, which again provides information about the nature of  the splittings of $F_n$ and hence of the complex of free factors related to these splittings. A relation between the simplicial automorphism group of the graph of free splittings of $F_n$ and $\Out(F_n)$ given by using sphere complex of $M$ can be found in \cite{souto1}.
\newline

Now let us fix a maximal sphere system $\Sigma$. This is a collection of isotopy classes of  disjoint and non-trivial 2-spheres in $M$ such that no two of which are isotopic and the complementary components are 3-punctured 3-spheres. By the intersection number of a torus $t$ we will mean the number $i(t,\Sigma)$ of components of intersection of $t$ with spheres of $\Sigma$ when the intersection is transversal. From now we assume that the intersections with spheres of $\Sigma$ are all transversal.

\begin{definition} We will say that a torus $t$ is  essential in $M$ if the image of  $\pi_1(t)$ is nontrivial in $\pi_1(M)$.
\end{definition}

The main goal of this paper is to define a ``normal form'' for an arbitrary homotopy class of an imbedded essential torus in $M$ with respect to a given sphere system,  which will amount to finding a representative with minimal intersection in the homotopy class. Naturally, we first need to show that such a nice representative exists in a given homotopy class of essential tori. This will be achieved by the following theorem:

\begin{theorem}\label{thm:main1}
Every imbedded essential torus in $M$ is homotopic to a normal torus and the homotopy process does not increase the intersection number with any sphere of $\Sigma$.
\end{theorem}
Secondly, we will prove the following theorem, which is analogous to \cite[Proposition 1.2]{H1} of Hatcher:
\begin{theorem}\label{thm:main2}
If $t$ and $t^\prime$ are two homotopic tori in $M$, both in normal form, then they are normally homotopic.
\end{theorem}
Here we define the notion of being \textit{normally homotopic} as follows:
\begin{definition} Two tori are said to be normally homotopic if there is a homotopy of $M$ changing one of the tori to the other one without introducing new intersections on the sphere crossings, hence through normal, but possibly immersed tori at each level.
\end{definition}

Using these two theorems, we will deduce:

\begin{corollary}\label{cor:1}
If a torus $t$ is in normal form with respect to a maximal sphere system $\Sigma$, then the intersection number of $t$ with any $S$ in $\Sigma$ is minimal among the representatives of the homotopy class $[t]$ in each $P$.
\end{corollary}

\subsection*{Motivation}
Hatcher in the paper \cite{H1} defined the notion of normal form with respect to a fixed sphere system and proved the existence of normal representatives of spheres in a given isotopy class of spheres in $M$. This leads to arguments about intersection numbers and minimal intersection conditions of these spheres and a correspondence between the free splittings of the free group $F_n$ and the embedded spheres in $M$ where $M=S^2\times S^1$, as in \cite{Gad2}.
\newline

We will extend this idea in this paper and look at imbedded  essential tori in $M$. This gives us a geometric interpretation of intersections of certain non-free group splittings and one might hope that it could be extended to other group splittings.
\newline

By a \textit{free splitting} of a group $G$, we mean the free product of two groups $A$ and $B$ so that $G$$=$$A\ast B$. By Kneser's conjecture, each splitting of the fundamental group of a closed three manifold as a free product corresponds to a sphere in the manifold. The details of this proof could be found in \cite{Hempel}.
\newline

By an \textit{amalgamated free product} of two groups $A$ and $B$ amalgamated along a group $C$, we mean the pushout of $A$ and $B$ when the maps $\alpha_1\colon C\rightarrow$$A$ and $\alpha_2\colon C\rightarrow$$B$  are group homomorphisms. It is denoted by $A\ast_C B$. When $C$ is a trivial group, we have a free product of $A$ and $B$. Similarly, an amalgamated free product of $A$ along $C$ ( or \textit{HNN extension} of $A$ along $C$) is defined to be the pushout of $A$ where $\iota_1,_2\colon C\rightarrow$$A$ are both homomorphisms and this universal group is denoted by $A\ast_C$. By a \textit{$\mathbb{Z}$--splitting} of a group $G$ we mean an amalgamated free product or an HNN extension of the group $G$ so that $C$ is isomorphic to $\mathbb{Z}$.
\newline

 It is known that the homotopy classes of essential imbedded tori in $M$ are in 1-1 correspondence with the equivalence classes of  $\mathbb{Z}$--splittings of $F_n$, namely, if the class of tori is separating it corresponds to an amalgamated free product and if it is nonseparating, it corresponds to an HNN extension of $F_n$. A proof of this is due to Matt Clay and can be found in \cite{CR}. The proof uses the connection between amalgamated free products (or HNN extensions) and the free splittings of $F_n$. These latter results are due to Shenitzer (amalgamated case) and Swarup (HNN case) and proofs can be found in \cite{STA}. This correspondence has some nice consequences which we consider as applications for our work. For example, there is an algebraically defined ``Dehn Twist'' notion on $\mathbb{Z}$--splittings of $F_n$ given in \cite{CLAYP}, which is used to generate fully irreducible elements of $\Out(F_n)$. This notion can be carried to tori in $M$ using this correspondence to find a topological way of generating fully irreducible elements, as in \cite{CR}. But before we attempt investigating such possible consequences, we need a notion of minimal intersection with respect to a given sphere system in $M$.
\newline

\subsection*{Acknowledgements}
I would like to thank Darryl McCullough for the suggestions, corrections and long discussions about the material throughout the preparation process which improved the text greatly and Kasra Rafi for giving the idea and the motivation for the paper and for his many useful suggestions.
\section{background}
We start by describing $M$ explicitly.
Since $M$ is a reducible 3-manifold, it has essential (not bounding 3-balls in $M$) 2-spheres which provide us with a rich algebraic structure. One way to describe $M$ is to remove the interiors of  $2n$ disjoint 3-balls from $S^3$ and identify the resulting 2-sphere boundary components in pairs by orientation-reversing diffeomorphisms, creating $S^2\times S^1$ summands.
 \newline

To give the connection between the two spaces on which $\Out(F_n)$ acts, first we need to define a few concepts. A \textit{sphere system}  is a collection of isotopy classes of disjoint and non-trivial  2-spheres in $M$ no two of which are isotopic. There is a simplicial complex associated to $M$ called the \textit{sphere complex} and denoted by $\SP(M)$, having isotopy classes of non-trivial 2-spheres in $M$ as vertices and sphere systems of $k+1$ spheres as $k$-dimensional simplices.
\newline

 The sphere complex has a subspace which is homeomorphic to Outer Space. This is the complement of the subcomplex which consists of sphere systems with at least one non simply-connected complementary component. For the details we refer to \cite{H1} and \cite{HV1}.
 \newline

 Throughout the paper, we will call 3-punctured 3-spheres either 3-punctured spheres or twice-punctured  3-cells(balls). These are analogous to pants in dimension 2. Hence we will denote them by $P$.
 \newline

 Using the essential spheres in $M$, we will be able to  define 3-dimensional versions of maximal curve systems and of pants decompositions for $M$ as follows:

\begin{definition}Given $M$ as above, we call a collection $\Sigma$ of disjointly imbedded essential, non-isotopic 2-spheres in $M$
a maximal sphere system if every  complementary component of $\Sigma$ in $M$ is a 3-punctured sphere.
\end{definition}

\section{tori in normal form}
Given an imbedded torus in $M$ and a maximal sphere system $\Sigma$, we can look at the number of intersections of the torus with the 2-spheres in each $P$, and define a notion of minimal intersection. In this paper we are particularly interested in the existence of a torus in a homotopy class which intersects the 2-spheres of a maximal sphere system $\Sigma$ minimally. There are certain pieces of a given torus in a $P$ that are particularly important for minimal intersection. They are:
\begin{enumerate}
\item  A disk piece, which is essential, in other words not parallel into any of the boundary 2-spheres and which has a single circle intersection with a single boundary sphere.
\item A cylinder piece, which is the topological boundary of a regular neighborhood of an arc connecting two different boundary 2-spheres.
\item A pants piece, which is the topological boundary of a regular neighborhood of a letter Y intersecting all three boundary components.
\end{enumerate}
\begin{definition}
Given an imbedded torus and a maximal sphere system $\Sigma$ in $M$, we say that the torus is in normal form with respect to $\Sigma$ if each intersection of the torus with each complementary 3-punctured sphere $P$ is a disk, a cylinder or a pants piece.
\end{definition}

Following Hatcher, we will show with the next theorem that any homotopy class of essential tori has a normal representative:

\begin{maintheorem1}
Every imbedded essential torus in $M$ is homotopic to a normal torus and the homotopy process does not increase the intersection number with any sphere of $\Sigma$.
\end{maintheorem1}
\begin{proof} The proof is analogous to the proof of Proposition 1.1 in \cite{H1} except that the homotopies are not necessarily isotopies as we do not have an analogue of Laudenbach's result for spheres given in \cite{LAUD} for tori, and consequently  we restrict ourselves to the homotopy classes instead of isotopy classes. As this result is of such importance in our work, we sketch the argument here.
 \newline

 Let us pick a representative $t$ from a homotopy class of tori.
 \newline

 As the first step, in each $P$, we regard each piece of the torus as consisting of sphere pieces inside $P$ and possibly concentric tubes connecting these sphere pieces to the boundary spheres of $P$. To do this, we first surger each piece of the torus along the intersection circles on the boundaries, starting from the innermost one, ending with the outermost one, resulting in a 2-sphere in $P$. On these sphere pieces, we reverse this surgery process by putting tubes between the sphere piece and the boundary spheres, in exactly the reverse order.
 \newline

  If $t$ is not normal, there will be a piece $F$ that meets a boundary sphere $S$ of a thrice punctured sphere $P$ in two intersection circles $C_1$ and $C_2$. Choose an arc $\alpha$ in $F$ connecting $C_1$ to $C_2$. Let $\alpha^\prime$ be an arc on $S\in$$\Sigma$ connecting $C_1$ to $C_2$. Reselecting $\alpha$ and $F$ if necessary, we may assume that interior of $\alpha^\prime$ does not meet $t$.
  \newline

   Let us call the tube portions of $F$ which meet $C_1$ and $C_2$ $T_1$ and $T_2$, respectively. Let us assume that $C_2$ was surgered before $C_1$.  Now, there is a homotopy of $t$ which is an isotopy on $F$ and whose effect is to slide the end of $T_1$ attached to the sphere part of $F$ along $\alpha$ to $T_2$ and finally out of $P$. Any tubes of $t$ inside of $T_1$ are slid along with it. If there are $r$ such tubes, at a certain point of the homotopy, they create $2r+1$ new intersection circles with $S$, $2r$ from the tubes inside $T_1$ and one from the intersection circle of the tube $T_1$ itself.
   \newline

   The intersection of $T_1$ with $P$ is now a cylinder. Since $\alpha$ is homotopic to $\alpha^\prime$, this cylinder along with any tubes inside it are homotopic to an imbedded position outside of $P$ near $\alpha^\prime$. During the homotopy, self intersections of $t$ may occur, but since the interior of $\alpha^\prime$ does not meet $t$, the final position of $t$ can be an imbedding. This homotopy eliminates $2r+2$ circles of intersection, giving a net decrease of 1 from the original position.
   \newline

 A sequence of such homotopies in each $P$ will give the desired homotopy in $M$. Since $t$ is essential, its image under a composition of such homotopies will not be disjoint form $\Sigma$ hence we will have a normal representative in the same homotopy class.
\end{proof}
\section{the universal cover}

A fixed maximal sphere system in $M$ gives a description of the universal cover $\widetilde{M}$ of $M$ as follows. Let $\PP$ be the set of twice punctured 3-balls in $M$ given by a maximal sphere system $\Sigma$ and regard $M$ as obtained from copies of $P$ in $\PP$ by identifying pairs of boundary spheres. Note that a pair might both be contained in a single $P$, in which case the image of $P$ in $M$ is a once-punctured $S^2\times S^1$. To construct $\widetilde{M}$, begin with a single copy of $P$ and attach copies of the $P$ in $\PP$ inductively along boundary spheres, as determined by unique path lifting. Repeating this process gives a description of $\widetilde{M}$ as a treelike union of copies of the $P$. We remark that $\widetilde{M}$ is homeomorphic to the complement of a Cantor set in $S^3$.
\newline

The universal cover $\widetilde{M}$ is modeled by a tree
$\widetilde{\Gamma}$, called the \textit{dual tree,} as follows. For each
copy $P$ in $\PP$ there is a vertex corresponding to the interior of $P$
and a vertex for each of the three boundary spheres, and there are three
edges connecting the interior vertex to the boundary sphere vertices. Hence there are two types of vertices: the valence-3 vertices indicating the 3-punctured spheres and valence-2 vertices indicating boundary spheres. To
obtain $\widetilde{\Gamma}$, identify the boundary vertices according to
how the corresponding sphere boundary components of the copies of the $P$
are identified to form $\widetilde{M}$.
\newline

We will call the union of the three edges for a copy of a $P$ a ``Y'',
since it is homeomorphic to a letter Y. We also write $\widetilde{\Sigma}$
for the union in $\widetilde{M}$ of the inverse images of the spheres in the
fixed sphere system.
\newline

Given a lift $\widetilde{t}$ of an imbedded torus in normal form, there is
a corresponding dual subgraph of $\widetilde{\Gamma}$ obtained by taking
the union of the Y's for the copies of the $P$ that meet
$\widetilde{t}$. We call this graph $T(\widetilde{t})$. The inclusion of
$\widetilde{t}$ into $\widetilde{M}$ is modeled by the inclusion of
$T(\widetilde{t})$ into $\widetilde{\Gamma}$, which is injective, hence we will have at most one component of $\widetilde{t}$ in each $P$. Note that
\begin{enumerate}
\item An extremal Y of $T(\widetilde{t})$, that
is, a Y that meets the rest of $T(\widetilde{t})$ in a single vertex,
occurs exactly when an intersection of $\widetilde{t}$ with a copy of a
$P$ is a disk. We will call such Y's type-1.
\item A Y meeting the rest of $T(\widetilde{t})$ in exactly two vertices
occurs exactly when an intersection of $\widetilde{t}$ with a copy of a
$P$ is a cylinder. These Y's will be called type-2.
\item A Y meeting the rest of $T(\widetilde{t})$ in its three boundary
vertices occurs exactly when an intersection of $\widetilde{t}$ with a
copy of a $P$ is a pair of pants. These are type-3.
\end{enumerate}
Since $\widetilde{t}$ is connected, $T(\widetilde{t})$ is also connected
and hence is a tree.
\section{the decorated graph}

To prove Theorem \ref{thm:main1}, we will provide a combinatorial description of the lift of a torus in terms of a tree defined in the universal cover of $M$. Such a lift equipped with a transverse orientation will be associated to a \textit{decorated graph}.
\newline

Let $t$ be an imbedded essential torus in $M$.
The image of $\pi_1(t)$ under the homomorphism induced by the inclusion $i\colon t\hookrightarrow$$M$ is an infinite cyclic subgroup of $\pi_1 (M)$, defined up to conjugacy. Fixing a specific lift of the inclusion to the universal cover, with image  $\widetilde{t}$, determines a specific subgroup in this conjugacy class, and a generator $\gamma$ of this subgroup acts as a covering transformation of $\widetilde{M}$ that preserves $\widetilde{t}$. Note that $\gamma$ does not interchange the sides of $\widetilde{t}$ since the image of $t$ is two-sided in $M$. There is a corresponding action of  $\pi_1(M)$ on $\widetilde{\Gamma}$ as simplicial isomorphisms. The generator $\gamma$ has an invariant axis which is topologically a line and $T(\widetilde{t})$ consists of this axis and finite trees meeting the axis. The action of $\gamma$ on $T(\widetilde{t})$ takes vertices to vertices and edges to edges. A fundamental domain for the action of $\gamma$ on $T(\widetilde{t})$ could be described as an arc on the invariant axis of $\gamma$ whose endpoints are translates by $\gamma$, together with some finite trees attached to this arc. Translates of these finite trees are all of the finite trees meeting the invariant axis of $\gamma$.
\newline

Motivated by Hatcher \cite{H1}, we will obtain the decorated graph with respect to the transverse orientation chosen as follows:
\newline

\begin{figure}
\setlength{\unitlength}{0.01\linewidth}
\begin{picture}(100,90)
\put(20,0){\includegraphics[width=0.7\textwidth]{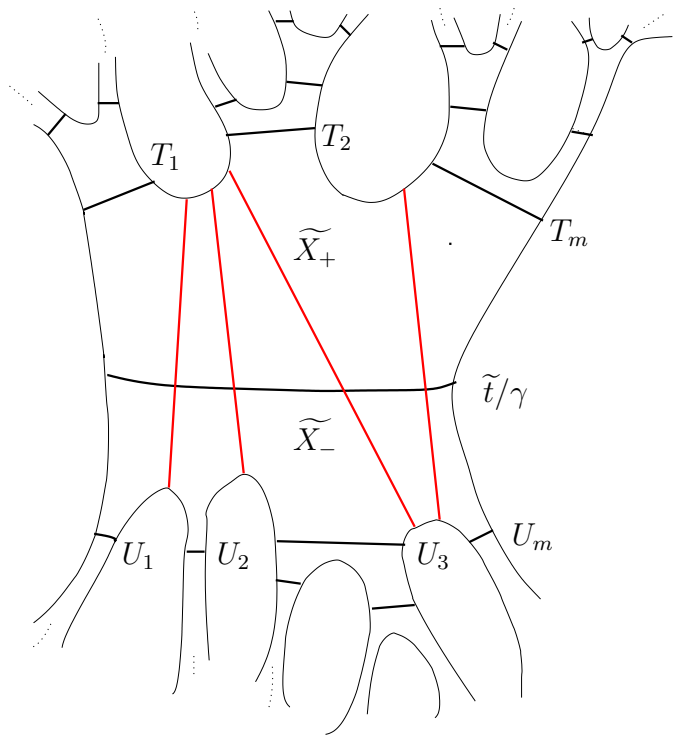}}
\put(42,18){$U_2$}
\put(32,18){$U_1$}
\put(63,18){$U_3$}
\put(73,20){$U_m$}
\put(35,60){$T_1$}
\put(53,62){$T_2$}
\put(77,52){$T_m$}

\put(50,50){$\widetilde{X_+}$}
\put(50,30){$\widetilde{X_-}$}
\put(70,35){$\widetilde{t}/\gamma$}

\end{picture}
\caption{A schematic picture of $\widetilde{t}/\gamma$, the submanifolds $\widetilde{X_+}$, $\widetilde{X_-}$, and spheres $T_i$$\in$$S_+$ and $U_j$$\in$$S_-$.}
\label{fig:HOM2}
\end{figure}

 We pick a transverse orientation of the lift $\widetilde{t}$ and label the sides; one with $+$ and the other $-$. This induces a corresponding orientation on $\widetilde{t}/\gamma$. Split $\widetilde{M}/\gamma$ along $(\widetilde{t}/\gamma)$$\cup$$(\bigcup$$\widetilde{S_i})$ where $\widetilde{S_i}$ are the spheres corresponding to the $\gamma$-orbits of the spheres in $\widetilde{M}$ which are disjoint from $\widetilde{t}$. Now, let $\widetilde{X_+}$ and $\widetilde{X_-}$ be the two components that contain copies of $\widetilde{t}/\gamma$ and define $S_+$=$\partial{\widetilde{X_+}} -\widetilde{t}/\gamma$ and $S_-$=$\partial{\widetilde{X_-}}-\widetilde{t}/\gamma$ where $\partial{\widetilde{X_+}}$ denotes the boundary of $\widetilde{X_+}$, etc. Note that $\widetilde{X_+}$,$\widetilde{X_-}$ and $\widetilde{t}/\gamma$ are compact submanifolds. See Figure \ref{fig:HOM2}.
 \newline

 We label the spheres $S_+$ with $+$ and the spheres $S_-$ with $-$. This gives a labeling of the vertices representing these spheres, which are extremal vertices of  $T(\widetilde{t})/\gamma$. For a disk piece of $\widetilde{t}/\gamma$, the corresponding two extremal vertices of $T(\widetilde{t})/\gamma$ will have the opposite signs. There will be no signs on a Y corresponding to a pants piece, since all spheres are intersected. For a  cylinder piece, one of the boundary spheres will not be intersected hence will be on one side of the torus and will correspond to a labeled extremal vertex on a type-2 Y. In particular if we have a torus with one of $S_+$ and $S_-$ empty, the corresponding graph will be some union of type-2 Y's, in other words finitely many extremal edges attached to the axis of $\gamma$ with one label on each of them, all labels the same. In this case, $t$ bounds a solid torus in $M$, $\widetilde{t}/\gamma$ bounds a solid torus $\widetilde{X_+}$ or $\widetilde{X_-}$ in $\widetilde{M}/\gamma$, and $\widetilde{t}/\gamma$ represents the trivial element in $H_2(\widetilde{M}/\gamma)$.
 \newline

\begin{figure}
\setlength{\unitlength}{0.01\linewidth}
\begin{picture}(100,60)
\put(20,0){\includegraphics[width=0.7\textwidth]{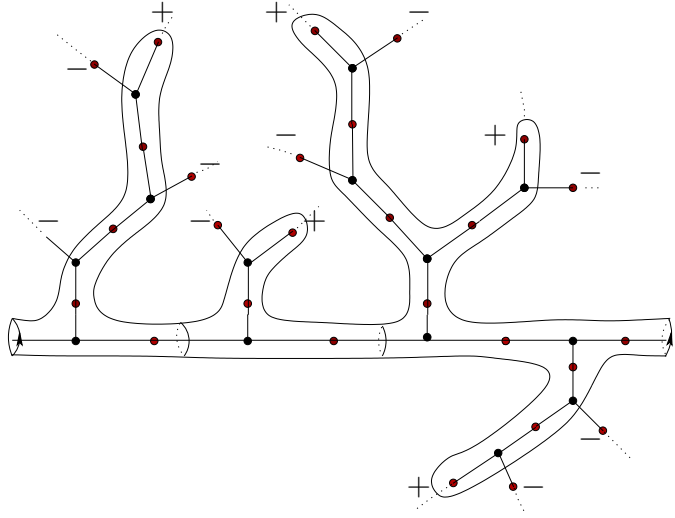}}
\put(51,30){$+$}
\put(47,52){$+$}
\put(26,46){$-$}
\put(40,36){$-$}
\put(23,30){$-$}
\put(35,52){$+$}
\put(39,30){$-$}
\put(48,39){$-$}
\put(70,39){$+$}
\put(80,35){$-$}
\put(80,7){$-$}
\put(74,2){$-$}
\put(62,2){$+$}
\put(62,52){$-$}

\end{picture}
\caption{The relation between $\widetilde{t}/\gamma$ and $g_t$.}
\label{fig:homb}
\end{figure}

The above construction  will give the graph $T(\widetilde{t})/\gamma$ in  $\widetilde{\Gamma}/\gamma$ a ``decoration'' of signs  on the ending vertices resulting from the transverse orientation on the torus $\widetilde{t}/\gamma$. We will call this decorated graph $g_t$ since it will be shown that it is uniquely determined by the normal homotopy class of the normal torus $t$. See Figure \ref{fig:homb}.

\section{the uniqueness of normal form}
Let $\Sigma$ be a maximal sphere system and $t$ and $t^\prime$ be two essential tori in $M$. For the following two lemmas, assume that the images of $\pi_1(t)$ and $\pi_1(t^\prime)$ in $\pi_1(M)$ are conjugate, say to the subgroup generated by $\gamma$.

\begin{lemma}\label{L1}
 Let $t$ and $t^\prime$ be two homotopic tori, both in normal form with respect to $\Sigma$. Then their transverse orientations may be chosen so that the corresponding decorated graphs with respect to the axis of $\gamma$ are equal.
\end{lemma}
\begin{proof} Let $t$ and $t^\prime$ be as given.
Define, as before, $S_+$=$\partial{\widetilde{X_+}}-\widetilde{t}/\gamma$, $S_-$=$\partial{\widetilde{X_-}}-\widetilde{t}/\gamma$, and $S^\prime_+$=$\partial{\widetilde{X_+}}-\widetilde{t^\prime}/\gamma$, $S^\prime_-$=$\partial{\widetilde{X_-}}-\widetilde{t^\prime}/\gamma$. We pick transverse orientations on $\widetilde{t}/\gamma$ and  $\widetilde{t^\prime}/\gamma$. These transverse orientations determine transverse orientations on $S_+$, $S_-$, $S^\prime_+$ and $S^\prime_-$ in $\widetilde{M}/\gamma$, and hence $+$ and $-$ labeling of them so that $\widetilde{t}/\gamma$ is homologous to both $S_+$, $S_-$, and $\widetilde{t^\prime}/\gamma$ is homologous to both $S^\prime_+$ and $S^\prime_-$.
\newline

 Now, any homotopy from $t$ to $t^\prime$ lifts to a homotopy from $\widetilde{t}/\gamma$ to $\widetilde{t^\prime}/\gamma$. Therefore we may fix transverse orientations on $t$ and $t^\prime$ so that $\widetilde{t}/\gamma$ and $\widetilde{t^\prime}/\gamma$ represent the same element of $H_2(\widetilde{M}/\gamma;\mathbb{Z})$. Then, $S_+$, $S_-$, $S^\prime_+$ and $S^\prime_-$ all represent the same homology class.
 \newline

Assume  that $g_t$ $\neq$ $g_{t^\prime}$.
Suppose first that $T(\widetilde{t})/\gamma$ $\neq$ $T(\widetilde{t^\prime})/\gamma$. Then one of them, say $T(\widetilde{t})/\gamma$ contains an extremal Y, say $Y_0$, not in $T(\widetilde{t^\prime})/\gamma$.
\newline

Consider the valence-2 vertex of $Y_0$ which connects it to the rest of the graph. Let us call it $v$. Now, $v$ represents a 2-sphere, which is a component of the boundary of the 3-punctured sphere $\widetilde{P}$ associated to the middle valence-3 vertex of the $Y_0$. This 2-sphere separates $\widetilde{M}$ into two parts. One part contains exactly one of the spheres in $S_+$ and one sphere in $S_-$ and the other part contains all of the spheres of $S^\prime_+$ and $S^\prime_-$ and all but the one of the spheres of $S_+$ and $S_-$. We will call this latter part $\widetilde{M_0}$. But then, $S^\prime_+$ and $S^\prime_-$  represent zero in $H_2(\widetilde{M}/\gamma,\widetilde{M_0}/\gamma)$ and $S_+$ and $S_-$ do not. This contradicts the fact that $S_+$, $S_-$, $S^\prime_+$ and $S^\prime_-$ represent the same homology class in $H_2(\widetilde{M}/\gamma)$.
  \newline

Now we are reduced to the case that $\t$ and $\tprime$ have the same
topological graphs. We must prove that their orientations may be selected
so that the decorations are equal.
\newline

Suppose first that $[\t]$ and hence $[\tprime]$ are $0$ in
$H_2(\widetilde{M}/\gamma)$. Then each bounds a compact submanifold of
$\widetilde{M}/\gamma$, so the decorations each have either all plus signs
or all minus signs. If they agree, there is nothing to prove. If they are
opposite, we may reverse the orientation on one of them (not changing its
homology class, since the class is $0$) to make the signs all agree, and
again we are finished. So we may assume that $[\t]$ is nonzero in $H_2(\widetilde{M}/\gamma)$.
\newline

We have $\partial\widetilde{X_+}=T_1+\cdots+T_m-\t$ and $\partial\widetilde{
X_-}=\t-U_1-\cdots-U_n$ where $T_i$$\in$$S_+$, $U_j$$\in$$S_-$ and $m,n$$\in$$\mathbb{Z}$. Since $[\t]$ is nonzero, $m$ and $n$ are both at
least~$1$.
\newline

We have $H_2(\widetilde{X_+},\t)\cong \mathbb{Z}^{m-1}=\langle (T_1\rangle \oplus\cdots \oplus
\langle T_m\rangle)/(T_1+\cdots +T_m=0)$. This is a subgroup of
$H_2(\overline{\widetilde{M}/\gamma-\widetilde{M_1}},\t),$ where
$\widetilde{M_1}$ is the component of $\widetilde{M}/\gamma$ cut along $\t$
that contains $\widetilde{X_-}$.
\newline

In fact, we have
\[ H_2(\widetilde{X_+},\t)\subset
H_2(\overline{\widetilde{M}/\gamma-\widetilde{M_1}},\t)\cong
H_2(\widetilde{M}/\gamma,\widetilde{M_1})\ ,\]
the latter isomorphism by excision, and under
\[ H_2(\widetilde{M}/\gamma)\to
H_2(\widetilde{M}/\gamma,\widetilde{M_1})\ ,\]
the homology class $[\t]$ goes into the subgroup
$H_2(\widetilde{X_+},\t)$ and equals $T_1+T_2+\cdots+T_m=0$.
\newline

Now $[\tprime]=[T_{i_1}+\cdots+T_{i_r}+U_{j_1}+\cdots+U_{j_s}]$
corresponding to the extremal vertices of the graph that are
decorated with plus signs for $\tprime$. Under
\[ H_2(\widetilde{M}/\gamma)\to
H_2(\widetilde{M}/\gamma,\widetilde{M_1})\ ,\]
$[\tprime]$ goes to $[T_{i_1}+\cdots+T_{i_r}]$, and must equal $0$ since it
equals $[\t]$ in $H_2(\widetilde{M}/\gamma)$. Therefore it contains either
all or none of the $T_i$. That is, in the decoration for the graph obtained
from $\tprime$, either all the $T_i$ have plus signs or all have minus signs.
\newline

Applying the same argument to the minus side (with $\widetilde{X_-}$ in the role of
$\widetilde{X_+}$), we conclude that for the decoration obtained from $\tprime$, either
all the $U_i$ have plus signs or all have minus signs. That is, we have
\[ [\tprime]=[\epsilon T_1\cdots +\epsilon T_m + \delta U_1+\cdots +\delta
  U_n]\ ,\]
where $\epsilon,\delta\in\{0,1\}$.
\newline

Suppose that $\epsilon=\delta=0$ or $\epsilon=\delta=1$, that is, in the
decoration for $\tprime$ all extremal vertices have either plus signs or
minus signs. Then $\tprime$ bounds a compact submanifold of
$\widetilde{M}$, contradicting the fact that $[\t]$ is nonzero.
\newline

If $\epsilon = 1$ and $\delta=0$, then the decorations are the same and
there is nothing to prove.
\newline

In the remaining case, when $\epsilon = 0$ and $\delta=1$, we may reverse
the orientation on $\tprime$ to make the decorations equal, and the proof
is complete.
\newline

Therefore, since the decorations agree also, we have $g_t$ = $g_{t^\prime}$.
\end{proof}

\begin{lemma}\label{L2}
For two tori $t$ and $t^\prime$ normal with respect to $\Sigma$, suppose that the corresponding decorated graphs are the same. Then, $t$ and $t^\prime$ are normally homotopic.
\end{lemma}
\begin{proof}

We will construct a normal homotopy of $t^\prime/\gamma$ in $\widetilde{M}/\gamma$, moving it onto $t/\gamma$. This projects to a normal homotopy of $t^\prime$ onto $t$ in $M$.
\newline

 Since the decorated graphs are the same, both tori will have the  same type of piece in each $\widetilde{P}$. To describe the normal homotopy, we start at one of the endpoints of the arc giving the axis of $\gamma$.
\newline

 Let us call the first $\widetilde{P}$ on the axis of $\gamma$, $\widetilde{P_1}$. By normal isotopy, we may move $t^\prime/\gamma\cap\widetilde{P_1}$ onto $t/\gamma\cap\widetilde{P_1}$. On the next $\widetilde{P}$ along the axis, say $\widetilde{P_2}$, we may move $t^\prime/\gamma\cap\widetilde{P_2}$ onto $t/\gamma$$\cap$$\widetilde{P_2}$ without moving $t^\prime/\gamma\cap(\widetilde{P_1}\cap\widetilde{P_2})$. It may be necessary to move $t^\prime/\gamma$ on the other components of $\partial{\widetilde{P_2}}$ using a ``twist''. Figure \ref{fig:twist} illustrates such a twist.
 \newline

 \begin{figure}
\begin{center}
\includegraphics[width=0.8\textwidth]{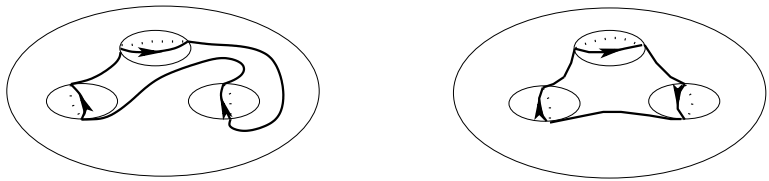}
\caption{A twist of one of the pieces on a boundary sphere.}
\label{fig:twist}
\end{center}
\end{figure}

 We continue along the axis of $\gamma$ in this way, until we reach $\widetilde{P_n}$ that meets $\widetilde{P_1}$. The isotopy moving $t^\prime/\gamma\cap\widetilde{P_n}$ onto $t/\gamma\cap\widetilde{P_n}$ can be accomplished without moving $t^\prime/\gamma\cap(\widetilde{P_n}\cap\widetilde{P_1})$, since if not, $t^\prime/\gamma$ would be a Klein Bottle.
 \newline

Now, we move to the finite tree branches on the axis of $\gamma$ and continue moving the pieces of $t^\prime/\gamma$  in each $\widetilde{P}$ corresponding to the Y's on branches, one after the other, fixing the already coinciding intersection circles we start with. Again, we might have the situation in Figure \ref{fig:twist}, so we might need to twist one intersection circle to make the pieces coincide. After a sequence of such homotopies we eventually reach an extremal Y, which must have one disk piece from each torus, one of the pieces with a twist, perhaps as in Figure \ref{fig:twist2}. Now, if we fill in the boundary spheres in this $\widetilde{P}$ with 3-cells, we will obtain a 3-ball, and by an isotopy we will be able to move one disk piece to the other one without moving boundary. Regarding the 3-cells as points, this determines an element of the braid group of two points in the 3-ball. This group is of order 2. But since the decorations are the same, the braid actually lies in the pure braid group of the 3-ball, which is trivial. So the disk pieces are isotopic relative to the boundary of $\widetilde{P}$.
\newline

\begin{figure}
\begin{center}
\includegraphics[width=0.8\textwidth]{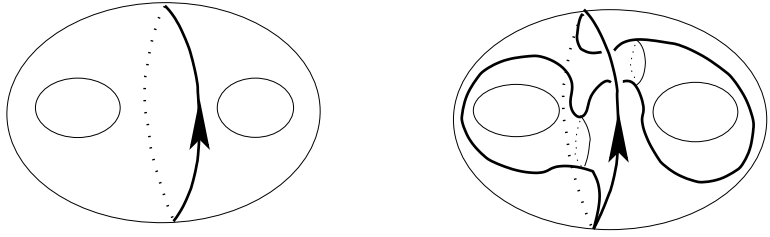}
\caption{A disk piece and a twisted disk piece.}
\label{fig:twist2}
\end{center}
\end{figure}

 We observe that at each stage of each of these isotopies, we have a normal torus because each piece only moves within a single $\widetilde{P}$.
As the last step, we take the composition of these isotopies and project it into $M$ to see that the two tori are normally homotopic. Since self intersections are possible, at some levels we might have immersed normal tori during this final homotopy.

\end{proof}
Now the proof of Theorem \ref{thm:main2}, which we will state here again, will be clear:
\begin{maintheorem2}
If $t$ and $t^\prime$ are two homotopic tori in $M$, both in normal form with respect to a maximal sphere system $\Sigma$, then they are normally homotopic.
\end{maintheorem2}
\begin{proof} Let $t$ and $t^\prime$ be two homotopic normal tori. We start with Lemma \ref{L1} to see that two tori have the same decorated graphs and continue with Lemma \ref{L2} to conclude that they are normally homotopic.
\end{proof}

\begin{corollary1}
If a torus $t$ is in normal form with respect to a maximal sphere system $\Sigma$, then the intersection number of $t$ with any $S$ in $\Sigma$ is minimal among the representatives of the homotopy class $[t]$ in each $P$.
\end{corollary1}
\begin{proof} Let $i(t,S)$ denote the number of components of $t$$\cap$$S$. Suppose $t$ is normal but there is a torus $t_1$ which is homotopic to $t$ with $i(t_1,S)< i(t,S)$. Then, by Theorem \ref{thm:main1}, $t_1$ is homotopic to a normal torus $t_2$ with $i(t_2,\Sigma)\leq i(t_1,\Sigma)$. Now, by Theorem \ref{thm:main2}
, any two homotopic normal tori are normally homotopic, which implies $i(t,S)$$=$ $i(t_2,S)\leq i(t_1,S)$$<i(t,S)$, a contradiction. Therefore  $i(t,S)$ was minimal among the tori in $[t]$.
\end{proof}

\bibliographystyle{alpha}
\bibliography{main}
\end{document}